\documentclass[11pt]{article}

\usepackage{amssymb, amsmath, amsthm, graphicx}
\usepackage[left=1in,top=1in,right=1in]{geometry}

\usepackage{subfigure}

      \theoremstyle{plain}
      \newtheorem{theorem}{Theorem}%[section]

      \newtheorem{proposition}[theorem]{Proposition}

      \theoremstyle{definition}

      \theoremstyle{remark}

\title{Coloring triple systems with local conditions}
\author{Dhruv Mubayi\thanks{Department of Mathematics, Statistics, and Computer Science, University of Illinois, Chicago, IL, 60607 USA.  Research partially supported by NSF grant DMS-1300138. Email: {\tt mubayi@uic.edu}} } 
\begin{document}
\maketitle
\medskip

\begin{abstract}
We produce an edge-coloring of the complete 3-uniform hypergraph on $n$ vertices with $e^{O(\sqrt{\log \log n})}$ colors such that the  edges spanned by  every set of five vertices receive at least three distinct colors. This answers the first open case of a question of Conlon-Fox-Lee-Sudakov~\cite{CFLS} who asked whether such a coloring exists with $(\log n)^{o(1)}$ colors.  
\end{abstract}

\section{Introduction}
A $k$-uniform hypergraph $H$ ($k$-graph for short) with vertex set $V(H)$ is a collection of $k$-element subsets of $V(H)$. Write $K_n^k$ for the complete $k$-graph with vertex set of size $n$.
A $(p,q)$-coloring of $K_n^k$ is an edge-coloring of $K_n^k$ that gives every copy of $K_p^k$ at least $q$ colors. Let $f_k(n,p,q)$ be the minimum number of colors in a $(p,q)$-coloring of $K_n^k$. This paper deals only with $k=3$. 

Conlon-Fox-Lee-Sudakov~\cite{CFLS} asked whether $f_3(n,p,p-2)=(\log n)^{o(1)}$ for $p \ge 3$ (the case $p=4$ is easy).
In this note we answer the first open case with a substantially smaller bound.

\begin{theorem}
$$f_3(n,5,3)= e^{O(\sqrt{\log\log n})}.$$
\end{theorem}
The problem of determining $f_k(n,p,q)$ for fixed $k,p,q$ has a long history, beginning with its introduction by Erd\H os and Shelah~\cite{E1, E2}, and  subsequent investigation (for graphs) by Erd\H os and Gy\'arf\'as~\cite{EG}. Studying $f_k(n,p,q)$ 
when $q=2$ is equivalent to studying classical Ramsey numbers, and most of the effort on these problems has therefore been for $q>2$. The simplest nontrivial  case in this regime is  $f_2(n,4,3)$, which was shown to be $n^{o(1)}$ in~\cite{M} and later $\Omega(\log n)$ (see~\cite{FS, KM}).  The same upper bound was shown for $f(n,5,4)$ in~\cite{EM}. Conlon-Fox-Lee-Sudakov~\cite{CFLS2} recently extended this construction considerably by proving  that
$f_2(n,p,p-1)=n^{o(1)}$ for all fixed $p \ge 4$. Their result is sharp in the sense that $f_2(n,p,p)=\Omega(n^{1/(p-2)})$. 

The first nontrivial hypergraph case is $f_3(n,4,3)$ and has tight connections to Shelah's breakthrough proof~\cite{S} of primitive recursive bounds for the Hales-Jewett numbers. Answering a question of Graham-Rothschild-Spencer~\cite{GRS}, Conlon et. al.~\cite{CFLS} recently proved that
 $f_3(n,4,3)= n^{o(1)}$. They also posed a variety of basic questions about $f_3(n,p,q)$,  including the one we address in this note. 
 
 Our construction  uses an extension of the coloring in~\cite{M} together with the stepping up technique of Erd\H os and Hajnal. It is quite possible that, similar to the situation for  graphs, other hypergraph cases will  eventually be addressed by the ideas introduced here. 

\section{The Construction} 

We begin by defining an edge-coloring $\sigma$ of the complete graph $K_n$ whose vertices are ordered.

{\bf Construction of $\sigma$:}  Given integers $t<m$ and $n={m \choose t}$, let $V(K_n)$ be the set of 0/1 vectors of length $m$ with exactly $t$ 1's. Write $v=(v(1), \ldots, v(m))$ for a vertex. The vertices are naturally ordered by the integer they represent in binary, so $v<w$ iff $v(i)=0$ and $w(i)=1$ where $i$ is the first position (minimum integer) in which $v$ and $w$ differ. By considering  vertices as characteristic vectors of sets, we may assume that 
$V(K_n)= {[m] \choose t}$ whenever convenient.  For each $B \in {[m]\choose t}$, let $f_B:2^B\rightarrow [2^t]$ be a bijection.  Given vectors $v<w$ that are characteristic vectors of sets $S<T$, let $c_1(vw)=\min \{i: v(i)=0, w(i)=1\}$, $c_2(vw)=\min \{j: j>i, v(i)=1, w(i)=0\}$, $c_3(vw)=f_S(S \cap T)$ and
$c_4(vw)=f_T(S \cap T)$. Finally, define $$\sigma(vw)=(c_1(vw), c_2(vw), c_3(vw), c_4(vw)).$$
If $n$ is not of the form ${m \choose t}$, then let $n'\ge n$ be the smallest integer of this form, color ${[n']\choose 2}$ as described above, and restrict the coloring to ${[n]\choose 2}$. \qed

It is known~\cite{M, M2} that $\sigma$ is both a $(3,2)$ and $(4,3)$-coloring of $K_n$ (we only need the first and fourth coordinates of color vectors  for this) and, for suitable choice of $m$ and $t$ it uses $e^{O(\sqrt{\log n})}$ colors for all $n$. We need the following additional properties.

\begin{proposition} The coloring $\sigma$ satisfies the following properties. 

1) If $v<w<x$, then $\sigma(vw)\ne \sigma(wx)$.

2) If $v<w<\min\{x,y\}$, and $\sigma(vw)=\sigma(vx)$, then $\sigma(vy)\ne \sigma(wx)$.

3) If $v<w<x<y$ with $\sigma(vw)=\sigma(xy)$, then 
$\sigma(vx)\ne \sigma(vy)$.
\end{proposition}

\proof It suffices to consider the first coordinate $c_1$ of $\sigma$ to prove the first two properties. For 1), observe that $i=c_1(vw)$ implies that $w(i)=1$, while $i=c_1(wx)$ implies that $w(i)=0$. For 2), let $i=c_1(vw)=c_1(vx)$ and suppose for contradiction that $i'=c_1(vy)=c_1(wx)$ so that $v(j)=y(j)$ for $j<i'$.  Assume first that $i<i'$  . Then $y(i)=v(i)=0$, while $w(i)=1$. This implies that $w>y$, a contradiction. Now assume that $i>i'$ ($i=i'$ is impossible since $w(i)=1$ while $w(i')=0$). Then  $0=v(i')=x(i')=1$ due to $c_1(vy)=i', c_1(vx)=i>i'$ and $c_1(wx)=i'$.

We now prove 3) so assume we are given  $v<w<x<y$ with  $c_1(vw)=c_1(xy)=i<j=c_2(vw)=c_2(xy)$.  Then $v(j)=x(j)=1$ and $y(j)=0$. Suppose that $v,w,x,y$ are characteristic vectors of $V, W, X, Y$ respectively. Then $c_3(vx)=c_3(VX)=f_V(V \cap X)$ while $c_3(vy)=c_3(VY)=f_V(V \cap Y)$. If  $c_3(vx)=c_3(vy)$, then $f_V(V \cap X)=f_V(V \cap Y)$ and since $f_V$ is a bijection, $V \cap X=V \cap Y$.  But this is impossible as $j \in (V \cap X)\setminus Y$.  \qed

We are now ready to describe the edge-coloring $\chi$ of $K_n^3$ that we will use. 

{\bf Construction of $\chi$:}  Given a copy of $K_n$ on $[n]$ and the edge-coloring $\sigma$, we produce an edge-coloring $\chi$ of the 3-graph $H$ on $\{0,1\}^n$ as follows. Order the vertices of $H$ according to the integer that they represent in binary. Given  vertices $x< y$ in $V(H)$, let $\gamma_{xy}$ be the first coordinate where $x$ and $y$ differ.  Given vertices $x<y<z$, let $\delta_{xyz}$ equal 1 if $\gamma_{xy}<\gamma_{yz}$ and $-1$ otherwise. For an edge $uvw$ with $u<v<w$, let
$$\chi(uvw)=(\sigma(\gamma_{uv} \gamma_{vw}), \delta_{uvw}). \quad \qed$$

Since $\sigma$ is an edge-coloring of $K_n$ with  $e^{O(\sqrt{\log n})}$ colors, $\chi$ is an edge-coloring of $K_N^3$ ($N=2^n$) with $e^{O(\sqrt{\log\log N})}$ colors as promised.  Moreover, extending this construction to all $N$ is trivial by considering the smallest $N'\ge N$ which is a power of 2, coloring ${[N']\choose 2}$ and restricting to ${[N]\choose 2}$. We are left with showing that $\chi$ is a $(5,3)$-coloring of $K_N^3$.

{\bf Proof that $\chi$ is a $(5,3)$-coloring:} Suppose, for contradiction, that $X=\{x_1, \ldots, x_5\}$ where $x_1<x_2<x_3<x_4<x_5$ are five vertices of $H$ forming a 2-colored $K_5^3$.  Let $\gamma_{i}=\gamma_{x_ix_{i+1}}$.  Let $\gamma=\min \gamma_j$ and assume this minimum is achieved by $\gamma_p$. Note that this minimum is uniquely achieved, and $\gamma_i\ne \gamma_{i+1}$ for all $i$.

{\bf Case 1:} $p\in \{1,4\}$. The arguments for both cases are almost identical so we only consider the case $p=1$.  By assumption we have $\gamma_1<\gamma_2$. First assume that $\gamma_3>\gamma_2$.  If $\gamma_4>\gamma_3$, then the $K_4$ on $\{\gamma_i: i \in [4]\}$ has three colors since $\sigma$ is a $(4,3)$-coloring and this gives at least three colors to the edges in $X$. If $\gamma_4<\gamma_3$ then the $K_3$ on $\{\gamma_i: i \in [3]\}$ has two colors since $\sigma$ is a $(3,2)$-coloring
and this gives  two colors to the edges of $H$ within $\{x_i: i \in [4]\}$ with positive $\delta$-coordinate. On the other hand $\delta_{x_3x_4x_5}=-1$, so we again have three colors on $X$. We now suppose that $\gamma_3<\gamma_2$.  If $\gamma_4<\gamma_3$, then 
the $K_3$ on $\{\gamma_2, \gamma_3, \gamma_4\}$ has two colors since $\sigma$ is a $(3,2)$-coloring
and this gives  two colors to the edges of $H$ within $\{x_2, x_3, x_4, x_5\}$ with negative $\delta$-coordinate. On the other hand $\delta_{x_1x_2x_3}=1$, so we again have three colors on $X$. Finally, we may assume that $\gamma_1<\gamma_3<\min\{\gamma_2, \gamma_4\}$. Now 
$\sigma(\gamma_1\gamma_3)\ne \sigma(\gamma_3, \gamma_4)$ due to property 1) of $\sigma$, hence
$\chi(x_1x_3x_4)\ne \chi(x_3x_4x_5)$ and both have positive $\delta$-coordinates. But $\delta_{x_2x_3x_4}=-1$, so $\chi(x_2x_3x_4)$  is the third color on $X$.

{\bf Case 2:}  $p\in\{2,3\}$. The arguments for both cases are almost identical so we only consider the case $p=2$. We have $\gamma_3> \gamma_2$.  If in addition $\gamma_4>\gamma_3$, then we get two colors among $\{x_2, x_3, x_4, x_5\}$ with positive $\delta$-coordinate while $\delta_{x_1x_2x_3}=-1$. So we may assume that $\gamma_2<\gamma_4<\gamma_3$. Now $\chi(x_2x_3x_4)$ and $\chi(x_2x_4x_5)$ both have positive $\delta$ coordinates while  $\delta_{x_3x_4x_5}=-1$. Hence we have three colors unless $\sigma(\gamma_2\gamma_3)=\sigma(\gamma_2\gamma_4)$ which we may assume. Certainly $\delta_{x_1x_2x_3}=-1$, so we are done unless $\sigma(\gamma_2\gamma_1)=\sigma(\gamma_4\gamma_3)$ which we also assume.  If $\gamma_1=\gamma_4$, then $\sigma(\gamma_2\gamma_4)=\sigma(\gamma_4\gamma_3)$ and hence $\{\gamma_2, \gamma_4, \gamma_3\}$ is a monochromatic triangle, contradiction. If $\gamma_1> \gamma_4$, then 
$\gamma_2<\gamma_4< \min\{\gamma_1,\gamma_3\}$ with 
$\sigma(\gamma_2\gamma_4)=\sigma(\gamma_2\gamma_3)$
and $\sigma(\gamma_2\gamma_1)=\sigma(\gamma_4\gamma_3)$.
This contradicts
property 2). 
If $\gamma_1< \gamma_4$, then 
$\gamma_2<\gamma_1<\gamma_4<\gamma_3$ with $\sigma(\gamma_2\gamma_1)=\sigma(\gamma_4\gamma_3)$ 
and $\sigma(\gamma_2\gamma_4)=\sigma(\gamma_2\gamma_3)$.
This contradicts property 3) and completes the proof.    \qed
\bigskip

\noindent \textbf{Acknowledgment.} I am grateful to David Conlon and Choongbum Lee for carefully reading an earlier draft of this note and giving comments that helped improve the presentation.

\end{document}